\documentclass{amsart}

\usepackage{amssymb}
\usepackage{amsmath}
\usepackage{graphicx}

\newtheorem{theorem}{Theorem}[section]

\newtheorem{proposition}[theorem]{Proposition}
\newtheorem{corollary}[theorem]{Corollary}

\newtheorem{_definition}[theorem]{Definition}
\newenvironment{definition}{\begin{_definition}\rm}{\end{_definition}}

\newtheorem{_remark}[theorem]{\it Remark}
\newenvironment{remark}{\begin{_remark}\rm}{\end{_remark}}

\newtheorem{_example}[theorem]{Example}

\numberwithin{equation}{section}
\numberwithin{table}{section}
\numberwithin{figure}{section}

\newcommand{\erase}[1]{}

\newcommand{\C}{\mathord{\mathbb C}}

\renewcommand{\P}{\mathord{\mathbb  P}}
\newcommand{\Q}{\mathord{\mathbb  Q}}
\newcommand{\R}{\mathord{\mathbb R}}
\newcommand{\Z}{\mathord{\mathbb Z}}

\newcommand{\BBB}{\mathord{\mathcal B}}

\newcommand{\OOO}{\mathord{\mathcal O}}

\newcommand{\QQQ}{\mathord{\mathcal Q}}

\newcommand{\TTT}{\mathord{\mathcal T}}

\newcommand{\ZZZ}{\mathord{\mathcal Z}}

\font\mathgot=eufm10

\newcommand{\ppp}{\mathord{\hbox{\mathgot p}}}

\newcommand{\maprightsp}[1]{\; \smash{\mathop{\longrightarrow}\limits\sp{#1}}\; }

\newcommand{\inj}{\hookrightarrow}
\newcommand{\surj}{\mathbin{\to \hskip -7pt \to}}
\newcommand{\isom}{\mathbin{\,\raise -.6pt\rlap{$\to$}\raise 3.5pt%
\hbox{\hskip .3pt$\mathord{\sim}$}\,}}
\newcommand{\mapdown}{\phantom{\Big\downarrow}\hskip -8pt \downarrow}

\newcommand{\set}[2]{\{\; {#1} \; \mid \; {#2} \;  \}}
\newcommand{\shortset}[2]{\{ {#1} \,|\, {#2}   \}}

\newcommand{\map}[3]{ #1 \, : \, #2 \, \to \, #3}
\newcommand{\shortmap}[3]{ #1 : #2 \to #3}

\newcommand{\sm}{\setminus}

\newcommand{\sprime}{\sp\prime}
\newcommand{\sptimes}{\sp{\times}}
\newcommand{\fiberproduct}{\Box}

\newcommand{\inv}{\sp{-1}}

\newcommand{\NS}{\mathord{\rm NS}}

\newcommand{\Ker}{\operatorname{\rm Ker}\nolimits}

\newcommand{\Gal}{\operatorname{\rm Gal}\nolimits}
\newcommand{\Sing}{\operatorname{\rm Sing}\nolimits}
\newcommand{\Spec}{\operatorname{\rm Spec}\nolimits}
\newcommand{\rank}{\operatorname{\rm rank}\nolimits}

\newcommand{\Emb}{\operatorname{\rm Emb}\nolimits}

\newcommand{\Tr}{\operatorname{\rm Tr}\nolimits}

\newcommand{\GL}{\operatorname{\it{GL}}\nolimits}
\newcommand{\SL}{\operatorname{\it{SL}}\nolimits}
\newcommand{\wt}[1]{\widetilde{#1}}

\newcommand{\ang}[1]{\langle #1\rangle}

\newcommand{\rmand}{\textrm{and}}
\newcommand{\quand}{\quad\rmand\quad}

\newcommand{\mystruth}[1]{\phantom{\hbox{\vrule height #1}}}
\newcommand{\mystrutd}[1]{\phantom{\hbox{\vrule depth #1}}}

\newcommand{\Cl}{\operatorname{\it Cl}\nolimits} 

\newcommand{\tllambda}{\tilde{\lambda}}

\newcommand{\tor}{{\mathord{\rm tor}}}
\newcommand{\tf}{{\mathord{\rm tf}}}
\newcommand{\Image}{\mathord{\rm Im}}
\newcommand{\XY}{{(X, Y)}}
\newcommand{\tub}{\TTT}

\begin{document}

\title[Non-homeomorphic conjugate varieties]
{Non-homeomorphic conjugate complex varieties}

\author{Ichiro Shimada}
\address{
Department of Mathematics,
Faculty of Science,
Hokkaido University,
Sapporo 060-0810,
JAPAN
}
\email{shimada@math.sci.hokudai.ac.jp
}


\begin{abstract}
We present a method to produce examples of
non-homeomorphic conjugate complex varieties
based on the genus theory of lattices.
As an application, we give examples of arithmetic Zariski pairs.
\end{abstract}

\maketitle

\section{Introduction}
We denote by $\Emb (\C)$ the set of embeddings $\sigma:\C\inj \C$
of the complex number field $\C$ into itself.
A \emph{complex variety}
is a reduced irreducible quasi-projective scheme over $\C$
with the classical topology.
For a complex variety $X$ and $\sigma\in \Emb(\C)$,
we define  a complex variety $X\sp\sigma$ 
by the following diagram of 
the fiber product:
$$
\renewcommand{\arraystretch}{1.2}
\begin{array}{ccc}
 X\sp\sigma & \maprightsp{} & X \\
 \mapdown &\fiberproduct & \mapdown\\
 \Spec\C &\maprightsp{\sigma^*} & \Spec\C.
\end{array}
$$
Two complex varieties $X$ and $X\sprime$
are said to be \emph{conjugate} if
there exists  $\sigma\in \Emb(\C)$
such that
$X\sp\sigma$ is isomorphic to  $X\sprime$ over $\C$.
It is easy to see that the 
relation of being conjugate is an equivalence relation.
\par
%
%
Serre~\cite{MR0166197}
gave an example of conjugate complex smooth projective varieties
that have different homotopy types.
After this famous result,
only few examples of non-homeomorphic conjugate complex  varieties seem to be  known
(e.g.~\cite{MR0349679}, \cite{MR2247887}).
The purpose of this note is to give a simple method
to produce  many explicit examples of non-homeomorphic conjugate complex varieties.
This method  is based on a topological idea 
in~\cite{AZP},
and the arithmetic theory of transcendental lattices
of singular $K3$ surfaces in~\cite{tssK3},
which has been generalized by~Sch\"utt~\cite{MS}.
\section{A topological invariant}
For a $\Z$-module $A$, we denote by $A_{\tor}$ the torsion part of $A$,
and by $A^{\tf}$ the torsion-free quotient $A/A_{\tor}$.
Note that 
a symmetric  bilinear form  $A\times A\to \Z$
naturally induces
a symmetric  bilinear form  $A\sp\tf\times A\sp\tf \to \Z$.
A \emph{lattice} is  a free $\Z$-module $L$
of finite rank  with a \emph{non-degenerate} symmetric bilinear form
$L\times L\to \Z$.
For a topological space $Z$,
let  $H_k(Z)$ denote the homology group $H_k(Z,\Z)$.
\par
\smallskip
Let $U$ be an oriented  topological manifold of dimension $4n$.
We denote by 
$$
\map{\iota_U}{H_{2n}(U)\times H_{2n}(U)}{\Z}
$$
the  intersection pairing.
We put
$$
J_\infty (U)\;:=\;\textstyle{\bigcap_K}\, \Image (H_{2n}(U\setminus K)\to H_{2n}(U)),
$$
where $K$ runs through the set of  compact subsets of $U$,
and $H_{2n}(U\setminus K)\to H_{2n}(U)$ is the natural homomorphism
induced by the inclusion.
We then put
$$
\wt{B}_U:=H_{2n}(U)/ J_\infty (U)
\quand
B_U:=(\wt{B}_U)^{\tf}.
$$
Since any topological cycle is compact,
the intersection pairing  $\iota_U$ induces  symmetric bilinear forms
$$
\map{{\tilde\beta}_U}{\wt{B}_U\times \wt{B}_U}{\Z} \quand \map{\beta_U}{B_U\times B_U}{\Z}.
$$
It is obvious that, if $U$ and $U\sprime$ are homeomorphic,
then there exists an isomorphism $(B_U, \beta_U)\cong (B_{U\sprime}, \beta_{U\sprime})$.
\par
\smallskip
Let $X$ be a smooth complex  projective variety of dimension $2n$.
Then
$H_{2n}(X)^\tf$ is a lattice by the intersection pairing $\iota_X$.
Let $Y_1, \dots, Y_m$ be irreducible subvarieties  of $X$ with codimension $n$.
We put
$$
Y:=Y_1\cup\dots\cup Y_m, \quad U:=X\sm Y, 
$$
and investigate the topological invariant $(B_U, \beta_U)$
of the complex variety $U$.
We denote by $\wt{\Sigma}_\XY$ the submodule of $H_{2n}(X)$ generated by 
the homology classes $[Y_i] \in H_{2n}(X)$, and put
 $\Sigma_\XY:=(\wt{\Sigma}_\XY)^\tf$.
 We then put
$$
\wt{\Lambda}_\XY:=\set{x\in H_{2n}(X)}{\iota_X(x, y)=0\;\;\textrm{for any}\;\; y\in \wt{\Sigma}_\XY},
$$
and  $\Lambda_{\XY}:=(\wt{\Lambda}_\XY)^{\tf}$.
Finally, we denote by
\begin{eqnarray*}
 \tilde{\sigma}_\XY : \wt{\Sigma}_\XY\times \wt{\Sigma}_\XY \;\to\; \Z,& \qquad& 
 \sigma_\XY : \Sigma_\XY\times \Sigma_\XY \;\to\; \Z, \\
 \tilde{\lambda}_\XY : \wt{\Lambda}_\XY\times \wt{\Lambda}_\XY \;\to\; \Z,
& \qquad& 
\lambda_\XY : \Lambda_\XY\times \Lambda_\XY \;\to\; \Z,
\end{eqnarray*}
the symmetric bilinear forms induced  by $\iota_X$.
\par
\smallskip
The proof of the following theorem  is essentially same as the argument in the proof of~\cite[Theorem 2.4]{AZP}.
We present a proof for the sake of completeness.
\begin{theorem}\label{thm:main}
Let $X$, $Y$ and $U$ be as above.
Suppose that $\sigma_\XY$ is non-degenerate.
Then $(B_{U}, \beta_{U})$ is isomorphic to $(\Lambda_\XY, \lambda_\XY)$.
\end{theorem}
\begin{proof}
Let $\tub\subset X$ be a tubular neighborhood of $Y$.
We  put
$\tub\sptimes :=\tub\sm Y=\tub\cap U$,
and denote by
\begin{eqnarray*}
& i_T\;:\; H_{2n}(\tub\sptimes)\to H_{2n}(\tub),  &  i_U\;:\; H_{2n}(\tub\sptimes)\to H_{2n}(U), \\
& j_T\;:\; H_{2n}(\tub)\to H_{2n}(X),\;\;  &  j_U\;:\; H_{2n}(U)\to H_{2n}(X), 
\end{eqnarray*}
the homomorphisms induced by the inclusions.
We first show that 
\begin{equation}\label{eq:ImjU}
\Image (j_U)=\wt{\Lambda}_\XY.
\end{equation}
It is obvious that  $\Image (j_U)\subseteq \wt{\Lambda}_\XY$.
Let $[W]\in \wt{\Lambda}_\XY$ be
represented by a real $2n$-dimensional topological cycle $W$.
We can assume that $W\cap Y$ consists of a finite number of points in $Y\sm \Sing (Y)$,
and that, locally around each intersection point $P$, the  topological cycle
$W$ is a differentiable manifold
intersecting $Y$ transversely at $P$.
Let $P_{i, 1}, \dots, P_{i, k(i)}$
(resp.~$Q_{i, 1}, \dots, Q_{i, l(i)}$)
be the intersection points of $W$ and $Y_i$ with local intersection number $1$ (resp.~$-1$).
Since $\iota_X([W], [Y_i])=0$, we have $k(i)=l(i)$.
For each $j=1, \dots, k(i)$, we choose a path
$$
\map{\xi_{i, j}}{I}{Y_i\sm \Sing(Y)}
$$
from $P_{i, j}$ to $Q_{i, j}$,
where $I:=[0, 1]\subset \R$ is the closed interval.
Let $\BBB$ denote a real $2n$-dimensional closed ball with the center $O$.
We can \emph{thicken} the path $\xi_{i, j}$ to a continuous map
$$
\map{\tilde{\xi}_{i, j}}{\BBB\times I}{X}
$$
in such a way that
$\tilde{\xi}_{i, j}\inv (Y)$ is equal to $\{O\}\times I$,
that the restriction of $\tilde{\xi}_{i, j}$ to $\{O\}\times I$ is equal to $\xi_{i, j}$,
and that the restriction of 
$\tilde{\xi}_{i, j}$ to $\BBB\times \{0\}$ (resp.~to $\BBB\times \{1\}$)
induces a homeomorphism from $\BBB$ to a closed neighborhood 
$\Delta\sp+_{i, j}$ of $P_{i, j}$ (resp.~$\Delta\sp-_{i, j}$ of $Q_{i, j}$) in $W$.
We then put
$$
W\sprime\;\;:=\;\; \bigl( \, W \,\sm \,\textstyle{\bigcup_{i, j}} \,(\Delta\sp+_{i, j} \cup \Delta\sp-_{i, j}) \bigr)\;\cup\; 
\textstyle{\bigcup_{i, j}}   \,\tilde{\xi}_{i, j}(\partial \BBB\times I).
$$
We can give an orientation to each $\tilde{\xi}_{i, j}( \BBB\times I)$
in such a way that $W\sprime$ becomes a topological cycle.
Then we have $[W]=[W\sprime]$ in $H_{2n}(X)$
and $W\sprime\cap Y=\emptyset$.
Therefore $[W]=[W\sprime]$  is contained in $\Image (j_U)$, and hence~\eqref{eq:ImjU} is proved.
Next we show that
\begin{equation}\label{eq:KerjU}
\Ker (j_U)\,\subseteq\, J_\infty (U).
\end{equation}
Consider the Mayer-Vietoris sequence 
$$
\cdots \;\maprightsp{}\;H_{2n}(\tub\sptimes)
\;\maprightsp{i}\;
H_{2n}(\tub)\oplus H_{2n}(U)
\;\maprightsp{j}\;
H_{2n}(X)\;\maprightsp{}\;\cdots, 
$$
where $i(x)=(i_T(x), i_U(x))$ and $j(y, z)=j_T(y)-j_U(z)$.
If $j_U(z)=0$, then $( 0, z)\in \Ker(j)=\Image (i)$,  and hence $z\in  \Image (i_U)$.
On the other hand, 
we have  $\Image (i_U)=J_\infty (U)$,
because $\tub$ is a tubular neighborhood of $Y$.
Hence~\eqref{eq:KerjU} is proved.
\par
\smallskip
Since $(H_{2n}(X)^\tf,\iota_X)$ is a lattice
and $(\Sigma_\XY, \sigma_\XY)$ is a sublattice by the assumption,
the orthogonal complement  $(\Lambda_\XY, \lambda_\XY)$ is also a lattice.
By~\eqref{eq:ImjU} and~\eqref{eq:KerjU}, we have a commutative diagram
\newcommand{\downsurj}{\raisebox{10pt}{\rotatebox{-90}{$\surj$}}}
\newcommand{\downinj}{\raisebox{8pt}{\rotatebox{-90}{$\inj$}}}
\newcommand{\verteq}{\raisebox{-2pt}{\rotatebox{90}{$=$}}}
\begin{equation}\label{eq:diag}
\renewcommand{\arraystretch}{1.32}
\begin{array}{cccccccccc}
0 & \longrightarrow  & \Ker (j_U) & \longrightarrow & H_{2n}(U) & \maprightsp{j_U} & \wt{\Lambda}_\XY & \longrightarrow & 0 &\textrm{(exact)\phantom{.}}\\
 & & \downinj  & & \verteq & & \phantom{v}\;\downsurj \;\;\lower 2pt \hbox{${}^{\tilde{v}}$} & \\
0 & \longrightarrow  &  J_\infty (U) & \longrightarrow & H_{2n}(U) & \longrightarrow & \wt{B}_U & \longrightarrow & 0 &\textrm{(exact).}\\
\end{array}
\end{equation}
By the definition of the intersection pairing, we have
$\iota_U(z, z\sprime)=\iota_X(j_U(z), j_U(z\sprime))$
for any $z, z\sprime \in H_{2n}(U)$.
Therefore the homomorphism $\tilde{v}$ in~\eqref{eq:diag} satisfies 
$$
\tilde{\lambda}_\XY (\zeta, \zeta\sprime)=\tilde{\beta}_U(\tilde{v}(\zeta),\tilde{v}(\zeta\sprime))
$$
for any $\zeta, \zeta\sprime \in \wt{\Lambda}_\XY$.
If $\tilde{v}(\zeta)\in (\wt{B}_U)_\tor$, then
$\zeta\in (\wt{\Lambda}_\XY)_\tor$ holds, because $\lambda_\XY$ is non-degenerate.
Hence $\tilde{v}\inv ((\wt{B}_U)_\tor)=(\wt{\Lambda}_\XY)_\tor$ holds.
Therefore $\tilde{v}$ induces an isomorphism $(\Lambda_\XY, \lambda_\XY)\cong (B_U, \beta_U)$.
\end{proof}
\section{Transcendental lattices}
A submodule $L\sprime$ of a free $\Z$-module $L$ is said to be \emph{primitive}
if $(L/L\sprime)_\tor=0$.
\par
\smallskip
Let $X$ be a smooth complex projective variety of
dimension $2n$.
Then we have a natural isomorphism $H_{2n} (X, \Z)^\tf \cong H^{2n}(X, \Z)^\tf$
that transforms  $\iota_X$ to the cup-product $(\phantom{i}, \phantom{\hskip 1pt})_X$.
Let  $S_X\subset H^{2n}(X, \Z)^\tf$ be
the submodule 
generated by the  classes $[Y]\in H^{2n}(X, \Z)^\tf$
of irreducible subvarieties $Y$ of $X$ with codimension $n$,
and let
$$ 
\map{s_X}{S_X\times S_X}{\Z}
$$
be the restriction of the cup-product to $S_X$.
Note that $s_X$ is non-degenerate
by Lefschetz decomposition  and Hodge-Riemann bilinear relations.
We consider the following  condition:
\begin{itemize}
\item[(P)]
$S_X$ is primitive in $ H^{2n}(X, \Z)^\tf$.
\end{itemize}
\begin{remark}
The condition (P) is satisfied if $\dim X=2$, because 
$S_X=H^2(X,\Z)^\tf\cap H^{1,1}(X)$ holds for a surface $X$.
For the case where $\dim X>2$,
see Atiyah-Hirzebruch~\cite{MR0145560} and Totaro~\cite{MR1423033}.
\end{remark}
Let $\sigma$ be an element of $\Emb(\C)$,
and consider the conjugate complex variety $X\sp\sigma$.
\begin{proposition}\label{prop:S}
The map $[Y]\mapsto [Y\sp\sigma]$
induces an isomorphism 
 $(S_X, s_X)\cong (S_{X\sp\sigma}, s_{X\sp\sigma})$.
\end{proposition}
\begin{proof}
Let $\ZZZ_X$ be the free $\Z$-module generated by irreducible subvarieties  $Y$ of codimension $n$ in $X$,
and let
$$
\map{\zeta_X}{\ZZZ_X\times\ZZZ_X}{\Z}
$$
be the intersection pairing.
Then $S_X$ is the image of the cycle map 
$Z\mapsto [Z]$ from $\ZZZ_X$ to $H^{2n}(X, \Z)\sp\tf$.
We put
$$
\BBB_X:=\set{Z\in \ZZZ_X}{\zeta_X(Z, W)=0\;\;\textrm{for any} \;\; W\in \ZZZ_X},
$$
and consider the \emph{numerical N\'eron-Severi lattice} $\NS_X:=\ZZZ_X/\BBB_X$ with the symmetric bilinear form
$\shortmap{\bar\zeta_X}{\NS_X\times\NS_X}{\Z}$
induced by $\zeta_X$.
Since $s_X$ is non-degenerate, 
the kernel of the cycle map $\ZZZ_X\to H^{2n}(X, \Z)\sp\tf$
coincides with $\BBB_X$,
and hence
$(S_X, s_X)$ is isomorphic to $(\NS_X, \bar\zeta_X)$.
In the same way,  we see that
$(S_{X\sp\sigma}, s_{X\sp\sigma})$ is isomorphic to $(\NS_{X\sp\sigma}, \bar\zeta_{X\sp\sigma})$.
On the other hand,
since the intersection pairing $\zeta_X$ is defined algebraically~(see Fulton~\cite{MR1644323}),
the map $Y\mapsto Y\sp\sigma$ induces an isomorphism 
$(\ZZZ_X, \zeta_X)\cong (\ZZZ_{X\sp\sigma}, \zeta_{X\sp\sigma})$, and hence 
it induces 
$(\NS_X, \bar\zeta_X)\cong (\NS_{X\sp\sigma}, \bar\zeta_{X\sp\sigma})$.
\end{proof}
\begin{definition}
We define the \emph{transcendental lattice} $T_X$ 
of $X$  by
$$
T_X:=\set{x\in H^{2n}(X, \Z)^\tf}{(x, y)_X=0\;\;\text{for any}\;\; y\in S_X}.
$$
\end{definition}
 \begin{theorem}\label{thm:main2}
Let $Y_1, \dots, Y_m$  be irreducible subvarieties of $X$
with codimension $n$
whose   classes $[Y_i]\in H^{2n}(X, \Q)$ span $S_X\otimes \Q$.
We put $Y:=\cup_{i=1}^m\,Y_i$ and $U:=X\sm Y$.
If $T_{X\sp{\sigma}}$ is not isomorphic to $T_{X}$,
then $U\sp{\sigma}$ is  not homeomorphic to $U$.
\end{theorem}
\begin{proof}
Note that
the classes  $[{Y_i}^\sigma]\in H^{2n}(X\sp\sigma, \Q)$  span $S_{X\sp\sigma}\otimes \Q$.
Theorem~\ref{thm:main} implies that
 $(B_U, \beta_U)$ is isomorphic to $T_X$, and 
$(B_{U\sp\sigma}, \beta_{U\sp\sigma})$ is isomorphic to $T_{X\sp\sigma}$.
Since $(B_U, \beta_U)$ is a topological invariant of $U$,
we obtain the result.
\end{proof}
\begin{proposition}\label{prop:genus}
Suppose that  {\rm (P)} holds  
for both of $X$ and $X\sp{\sigma}$.
Then 
$T_{X\sp\sigma}$ is contained in the same genus as $T_X$.
\end{proposition}
\begin{proof}
By the GAGA principle,
the Hodge numbers 
$h^{p, q}(X)=\dim H^q (X, \Omega^p)$
of a smooth projective variety $X$
are   invariant under the  conjugation.
Since the signature of the cup-product on $H^{2n}(X, \R)$
is given by these $h^{p, q}(X)$ (see, for example, \cite[Theorem 6.33]{MR1967689}),
the cup-products on $H^{2n}(X, \R)$ and on $H^{2n}(X\sp\sigma, \R)$
have the same signature.
\par
\medskip
By the condition (P) for $X$, 
we see that $S_X$ is a primitive sublattice of $H^{2n} (X, \Z)\sp\tf$. 
Since $H^{2n}(X, \Z)^\tf$ is unimodular,
the discriminant form of $T_X$ is isomorphic to 
the discriminant form of $S_X$ multiplied by $-1$
by Nikulin's result~\cite[Corollary~1.6.2]{MR525944}.
In the same way,
we see that the discriminant form of $T_{X\sp\sigma}$ is isomorphic to $-1$ times 
the discriminant form of $S_{X\sp\sigma}$.
Since $S_X$ and $S_{X\sp\sigma}$ are isomorphic
by~Proposition~\ref{prop:S},
the discriminant forms of $T_X$ and $T_{X\sp\sigma}$ are isomorphic.
Moreover, the signatures of $T_X$ and $T_{X\sp\sigma}$ are the same
by the argument in the previous paragraph.
Hence
$T_X$ and $T_{X\sp\sigma}$  are  contained in the same genus
by~\cite[Corollary~1.9.4]{MR525944}.
\end{proof}
\section{Arithmetic of singular abelian and $K3$  surfaces}
A lattice $L$ is called \emph{even} if $(v, v)\in 2\Z$ holds for any $v\in L$.
\par
\smallskip
Let  $X$ be  
a complex abelian surface (resp.~a complex algebraic $K3$ surface).
We say that $X$ is \emph{singular}
if  $\rank (S_X)$ attains the possible maximum $4$
(resp.~$20$).
Suppose that   $X$ is singular in this sense.
Then we have
$$
T_{X}=H^2(X, \Z)\cap (H^{2,0}(X)\oplus H^{0,2}(X)),
$$
and   $T_{X}$ is an even  positive-definite   lattice of rank $2$.
Moreover, $T_{X}$ has a canonical orientation given as follows;
an ordered basis $(e_1, e_2)$ of $T_{X}$ is \emph{positive}
if the imaginary part of  $(e_1, \omega)_X/(e_2, \omega)_X\in \C$
is positive, where $\omega$ is a basis of $H^{2, 0}(X)$.
We write   $\wt{T}_X$ for the \emph{oriented} transcendental  lattice of $X$.
\begin{theorem}[Shioda-Mitani~\cite{MR0382289}]\label{thm:ShiodaMitani}
The map $A\mapsto \wt{T}_A$ induces a bijection
from the set of isomorphism classes of complex singular abelian surfaces $A$
to the set of isomorphism classes of even positive-definite oriented-lattices of rank $2$.
\end{theorem}
For complex singular abelian surfaces,
we have a converse of Proposition~\ref{prop:genus}.
\begin{proposition}\label{prop:A}
Let $A$ be a complex singular abelian surface, and 
$\wt{T}\sprime$  an   even  positive-definite oriented-lattice of rank $2$
such that the underlying lattice $T\sprime$ is 
contained in the  genus of $T_A$.
Then there exists $\sigma\in \Emb(\C)$ such that $\wt{T}_{A\sp\sigma}\cong \wt{T}\sprime$.
\end{proposition}
In~\cite{tssK3}, we  proved Proposition~\ref{prop:A}
under  additional conditions.
Then Sch\"utt~\cite{MS} proved Proposition~\ref{prop:A}
in this full generality.
We present a proof that does not make use of the id\'ele groups.
\par
\medskip
For the proof, we fix notation and prepare theorems
in the theory of complex multiplications~(\cite{MR1028322}, \cite{MR890960}, \cite{MR0314766}).
For $a,b,c\in \Z$, we put 
$$
Q[a,b,c]:=
\left[
\begin{array}{cc}
 2 a & b \\
 b & 2c
\end{array}
\right].
$$
For a negative integer $d$,
we denote by $\QQQ_d$ the set of matrices $Q[a,b,c]$
such that $a,b,c\in \Z$, $a>0$, $c>0$ and $d=b^2-4ac$.
The group $\GL_2 (\Z)$ acts on $\QQQ_d$  by
$Q\mapsto {}^t g Q g$,
where $Q\in \QQQ_d$ and $g\in GL_2(\Z)$.
Then the set of orbits $\QQQ_d/ GL_2(\Z)$
(resp.~$\QQQ_d/SL_2(\Z)$)
is identified with
the set of isomorphism classes of even positive-definite lattices 
(resp.~oriented-lattices) of rank $2$ with discriminant $-d$.
For an $\SL_2(\Z)$-orbit $\Lambda\in \QQQ_{d}/\SL_2(\Z)$
and a positive integer $m$,
we put
$$
\ang{m} \Lambda := \set{Q[ma,mb,mc]}{Q[a,b,c]\in \Lambda}\;\;\in \;\;  \QQQ_{d m^2}/\SL_2(\Z).
$$
We denote by  $\QQQ\sp*_d$ the subset of $\QQQ_d$ consisting of matrices $Q[a,b,c]\in \QQQ_d$
with  $\gcd(a,b,c)=1$.
Then  $\QQQ\sp*_d$ is stable under the action of $\GL_2(\Z)$.
\par
\smallskip
Let $K\subset \C$ be an imaginary quadratic field.
We denote by $\Z_K$ the ring of integers of $K$,
and by $D_K$ the discriminant of $\Z_K$.
For a positive integer $f$,
let  $\OOO_f\subseteq \Z_K$ denote the order of conductor $f$.
By a \emph{grid},
we mean a $\Z$-submodule of $K$ with rank $2$.
For grids $L$ and $L\sprime$,
we write $[L]=[L\sprime]$ if $L=\lambda L\sprime$
holds for some $\lambda\in K\sp\times$.
We put
$$
\OOO(L):=\set{\lambda\in K}{\lambda L\subseteq L},
\quad f(L):=[\Z_K:\OOO(L)],
\quad d(L):= D_K f(L)^2.
$$
By  definition,  we have $\OOO(L)=\OOO_{f(L)}$.
Moreover, if $[L]=[L\sprime]$, then $f(L)=f(L\sprime)$.
We equip $L$ with a symmetric bilinear form  defined by 
$$
(x, y)_L:=\frac{n^2}{[\OOO_{f(L)}: n L]}\Tr_{K/\Q}(x\bar y),
$$
where $n$ is a positive integer such that $nL\subset \OOO_{f(L)}$.
It  turns out that $(\phantom{i}, \phantom{\hskip 1pt})_L$ 
takes values in $\Z$.
An ordered basis $[\alpha,\beta]$ of $L$ is defined to be \emph{positive}
if the imaginary part of $\alpha/\beta\in \C$ is positive.
With $(\phantom{i}, \phantom{\hskip 1pt})_L$ and this orientation,
  $L$ becomes  an even positive-definite oriented-lattice of discriminant $-d(L)$.
We denote by $\tllambda (L)\in \QQQ_{d(L)}/\SL_2(\Z)$
the isomorphism class of the oriented-lattice $L$,
and by $\lambda (L)\in \QQQ_{d(L)}/GL_2(\Z)$
the isomorphism class of the underlying lattice.
It is obvious that, if $[L]=[L\sprime]$,
then  $\tilde\lambda(L)=\tllambda (L\sprime)$ holds.
\par
\smallskip
Let $L$ and $M$ be grids.
We denote by $LM$ the grid generated by $xy$, 
where $x\in L$ and $y\in M$.
Then
$[L][M]:=[LM]$
is well-defined.
It is easy to prove that
\begin{equation*}\label{eq:fLM}
f(LM)=\gcd(f(L), f(M)).
\end{equation*}
Let $f$ be a positive integer.
We put $d:=D_K f^2$,
and consider the set
$$
\Cl_d:=\set{[L]}{f(L)=f},
$$
which is  a finite abelian group by the product $[L][M]=[LM]$
with the unit element $[\OOO_f]$ and the inversion $[L]\inv=[L\inv]$,
where $L\inv:=\shortset{\lambda\in K}{\lambda L\subseteq \OOO_f}$.
\begin{proposition}[\S3 and \S7 in \cite{MR1028322}]\label{prop:Cld}
{\rm (1)}  The map $L\mapsto \tllambda(L)$ induces a bijection $\Cl_d \cong \QQQ_d\sp * /\SL_2(\Z)$
with the inverse map  being induced from  $Q\mapsto [L_Q]$, where 
\begin{equation}\label{eq:LQ}
L_Q:= \Z +\Z\left(\frac{-b+\sqrt{d}}{2a}\right)
\;\;\textrm{for\;\; $Q:=Q[a,b,c]\in \QQQ_d\sp *$}.
\end{equation}

{\rm (2)} For grids $L$ and $M$ with $f(L)=f(M)=f$,
the lattices $\lambda (L)$ and $\lambda (M)$
are contained in the same genus if and only if $[L][M]\inv \in (\Cl_d)^2$ holds.
\end{proposition}
Let $I\subseteq\OOO_f$ be an $\OOO_f$-ideal.
Then $f(I)$ divides $f$.
We say that $I$ is a \emph{proper} $\OOO_f$-ideal if $f=f(I)$ holds.
For a non-zero integer $\mu$, 
we say that $I$ is \emph{prime to $\mu$} if $I+\mu\OOO_f=\OOO_f$ holds.
\begin{proposition}[Chapter~8 of \cite{MR890960}]\label{prop:JI}\label{prop:prime}
{\rm (1)}
Any $\OOO_f$-ideal prime to $f$ is proper.
The map $J\mapsto I=J\cap \OOO_f$ is a bijection from the set of
$\Z_K$-ideals $J$ prime to $f$ to the set of $\OOO_f$-ideals $I$  prime to $f$.
The inverse map is given by $I\mapsto J=I\Z_K$.

{\rm (2)}
Let $\mu $ be a non-zero integer.
For any $[M]\in Cl_d$,
there exists a proper $\OOO_f$-ideal $I$ prime to $\mu$
such that $[I]=[M]$.
\end{proposition}
For a grid $L$, we denote by $j(L)\in \C$ the $j$-invariant of 
the complex elliptic curve $\C/L$.
It is obvious that $j(L)=j(L\sprime)$ holds if and only if $[L]=[L\sprime]$.
We then put
$$
H_d:=K(j(\OOO_f)),\quad \textrm{where $d:=D_Kf^2$.}
$$
The set $\shortset{j(L)}{f(L)=f}$ is 
contained in $H_d$, and coincides with 
the set of conjugates of $j(\OOO_f)$ over $K$.
Hence $H_d/K$  is a finite Galois extension,
which is called the \emph{ring class field of $\OOO_f$}.
\begin{theorem}[Chapter~10 of \cite{MR890960}]\label{thm:CM}
{\rm (1)} 
We define  $\varphi_d: \Cl_d\to \Gal(H_d/K)$ by
$$ j(\OOO_f)^{\varphi_d ([M])}:= j(M\inv )
\;\;\textrm{for}\;\;[M]\in \Cl_d.
$$ 
Then  we have 
$ j(L)^{\varphi_d ([M])}= j(M\inv L)$
for any $[M], [L]\in \Cl_d$, and 
$\varphi_d$ induces  an isomorphism  $\Cl_d\cong \Gal(H_d/K)$.

{\rm (2)}
If a prime $\ppp\subset\Z_K$ of $K$ ramifies in $H_d$,
then $\ppp$ divides $f\Z_K$.
For a $\Z_K$-ideal  $J$ prime to $f$,  
 the Artin automorphism $(J, H_d/K)\in \Gal (H_d/K)$
 is equal to $\varphi_d ({[J\cap \OOO_f]})$.
\end{theorem}
Shioda and Mitani~\cite{MR0382289} proved the following.
Let $\wt{T}$ be an even positive-definite oriented-lattice of rank $2$ given  by 
$Q[a,b,c]\in \QQQ_d$.
We put 
$$
K:=\Q(\sqrt{d}), \;\; 
m:=\gcd(a,b,c),\;\; 
d_0:=d/m^2, \;\; 
Q_0:=Q[a/m,\,b/m,\,c/m]\in \QQQ_{d_0}\sp*.
$$
There exists a positive integer $f$ such that $d=D_K(mf)^2$.
We consider the grid $L_0:=L_{Q_0}$ of $K$
associated with $Q_0\in \QQQ_{d_0}\sp*$  by~\eqref{eq:LQ}.
Then we have $f(L_0)=f$ and $d(L_0)=d_0$.
Note that  
 $\ang{m}\tilde{\lambda}(L_0)\in \QQQ_{d}/\SL_2(\Z)$ is 
the isomorphism class  containing $\wt{T}$.
Note also that we have
$\OOO_{mf}=\Z+\Z\,({b+\sqrt{d}}\,)/{2}$.
For grids $L$ and $M$ of $K$, we denote by $A(L, M)$ the  complex abelian surface
$\C/L\times \C/M$. 
It is well-known that the elliptic curve $\C/L$
is defined over the subfield $\Q(j(L))$ of $\C$.
Hence $A(L, M)$
is defined over $\Q(j(L), j(M))$
\begin{theorem}[\cite{MR0382289}]\label{thm:SM2}
{\rm (1)}
The oriented transcendental lattice of $A(L_0, \OOO_{mf})$ is isomorphic to 
 $\wt{T}$.
In other words, $\wt{T}_{A(L_0, \OOO_{mf})}$ is contained in 
 $\ang{m}\tilde\lambda(L_0)$.
 
 {\rm (2)}
 Let $L_1$ and $L_2$ be grids
of $K$.
Then $A(L_1, L_2)$ is isomorphic to $A(L_0, \OOO_{mf})$
if and only if  $[L_1 L_2]=[L_0]$ and
 $f(L_1)f(L_2)=mf^2$ hold.
\end{theorem}
We are now ready to prove Proposition~\ref{prop:A}.
\begin{proof}[Proof of Proposition~\ref{prop:A}]
Since $T_A$ and $T\sprime$ are in  the same genus,
they have the same discriminant,
which we denote by  $-d$.
Let $\wt{T}_A$ and $\wt{T}\sprime$ be represented  by 
$$
Q_A=Q[a,b,c]\in \QQQ_d
\quand
Q\sprime=Q[a\sprime,b\sprime,c\sprime]\in \QQQ_d, 
$$
respectively.
Since $T_A$ and $T\sprime$ are in  the same genus,
we have
$$
m:=\gcd (a,b,c)=\gcd(a\sprime, b\sprime, c\sprime).
$$
As above, we put $K:=\Q(\sqrt{d})$,
$d_0:=d/m^2$, and let $f$ be the positive integer such that
$d=D_K(mf)^2$.
Let us consider the grids
$$
L_0:=L_{Q[a/m,\, b/m, \,c/m]}
\quand
L\sprime_0:=L_{Q[a\sprime\hskip -1pt /m, \,b\sprime\hskip -1pt /m, \,c\sprime\hskip -1pt /m]}
$$
associated with $(1/m)Q_A\in \QQQ\sp *_{d_0}$ and 
$(1/m)Q\sprime \in \QQQ\sp *_{d_0}$
by~\eqref{eq:LQ}.
We have 
$f(L_0)=f(L\sprime_0)=f$.
The oriented-lattices  $\wt{T}_A$ and  $\wt{T}\sprime$  are contained in the
isomorphism classes  $\ang{m}\tilde{\lambda} (L_0)\in \QQQ_d/\SL_2(\Z)$
 and $\ang{m}\tilde{\lambda} (L\sprime_0) \in \QQQ_d/\SL_2(\Z)$, respectively.

By Theorems~\ref{thm:ShiodaMitani} and~\ref{thm:SM2},
we have an isomorphism
$A\cong A(L_0, \OOO_{mf})$
over $\C$.
Note that $A(L_0, \OOO_{mf})$ is defined over $H_d \subset \C$,
because we have
$j(L_0)\in H_{d_0}\subseteq H_d$ by~\cite[Theorem 6 in Chapter 10]{MR890960}.

Since $T_A$ and $T\sprime$ are in  the same genus,
the lattices 
${\lambda} (L_0), {\lambda} (L\sprime_0)\in \QQQ\sp *_{d_0}/\GL_2(\Z)$ are also in the same genus.
By Propositions~\ref{prop:Cld}~and~\ref{prop:prime},  
there exists a proper $\OOO_f$-ideal $I_f$ prime to $mf$ 
such that $[L_0\sprime]=[I_f]^2[L_0]$ holds in $\Cl_{d_0}$.
We put $J:=I_f\Z_K$ and $I_{mf}:=J\cap\OOO_{mf}$.
Then $J$ is a $\Z_K$-ideal prime to $mf$
satisfying $J\cap \OOO_f=I_f$, and 
$I_{mf}$
 is a proper $\OOO_{mf}$-ideal prime to $mf$.
The Artin automorphism
$$
\tau:=(J, H_d/K)\;\in \;\Gal(H_d/K)
$$
is equal to $\varphi_d ([I_{mf}])$, and 
its restriction  to $H_{d_0}$ is equal to $\varphi_{d_0} ([I_{f}])\in \Gal(H_{d_0}/K)$.
We extend ${\tau\inv}\in \Gal(H_d/K)$ to $\sigma\in \Emb (\C)$.
Then  we have $j(L_0)\sp{\sigma}=j(I_{f}  L_0)$ and $j(\OOO_{mf})\sp{\sigma}=j(I_{mf})$.
Hence
$A(L_0, \OOO_{mf})\sp{\sigma}$
is isomorphic to 
$A(I_{f}   L_0, I_{mf} )$,
which is then isomorphic to 
$A(I_{mf}I_{f}  L_0, \OOO_{mf} )$
by Theorem~\ref{thm:SM2}.
We have 
\begin{equation}\label{eq:IIL}
[I_{mf} I_{f}  L_0]=[I_f]^2[L_0]=[L_0\sprime].
\end{equation}
To prove this, it is enough to show that
$I_{mf}\OOO_f=I_f$.
The inclusion $I_{mf}\OOO_f\subseteq I_f$ is obvious.
Since $I_{mf}$ is prime to $mf$,
 we have 
$$
I_f=I_f\OOO_{mf}=I_f(I_{mf} +mf  \OOO_{mf})\subseteq I_{mf} \OOO_f + m  I_f.
$$
Since $mI_f\subseteq m\OOO_f\subseteq \OOO_{mf}$,
we have $m  I_f\subseteq I_f\cap \OOO_{mf}=I_{mf}$.
Therefore $I_f \subseteq I_{mf} \OOO_f$ holds,
and~\eqref{eq:IIL} is proved.
Consequently,  the oriented transcendental lattice of 
$$
A\sp\sigma \;\cong\; A(L_0, \OOO_{mf})\sp{\sigma}\;\cong\; 
A(I_{mf}I_{f}   L_0, \OOO_{mf})\;\cong\; A(L_0\sprime, \OOO_{mf})
$$
is contained in $\ang{m}\tilde{\lambda} (L\sprime_0)$ by Theorem~\ref{thm:SM2},
and hence  is  isomorphic to $\wt{T}\sprime$.
\end{proof}
For  complex singular $K3$ surfaces, we have the following result:
\begin{theorem}[Shioda-Inose~\cite{MR0441982}]\label{thm:ShiodaInose}
The map $Y\mapsto \wt{T}_Y$ induces a bijection
from the set of isomorphism classes of complex singular $K3$ surfaces $Y$
to the set of isomorphism classes of even positive-definite oriented-lattices of rank $2$.
\end{theorem}
Let $Y$ be a complex singular $K3$ surface,
and  $A$ the complex singular abelian surface 
such that $\wt{T}_Y\cong \wt{T}_A$.
Then $Y$ is obtained from $A$ by Shioda-Inose-Kummer construction
(see~\cite[\S6]{tssK3}).
By~\cite[Proposition 6.4]{tssK3}, 
we have 
$\wt{T}_{Y\sp\sigma}\cong \wt{T}_{A\sp\sigma}$ for any $\sigma\in \Emb(\C)$.
Combining this fact with 
 Proposition~\ref{prop:A},
we obtain the following:
\begin{proposition}\label{prop:K3}
Let $Y$ be a complex singular $K3$  surface, and 
$\wt{T}\sprime$  an   even  positive-definite oriented-lattice of rank $2$
such that  $T\sprime$ and $T_Y$ are 
contained in the same genus.
Then there exists $\sigma\in \Emb(\C)$ such that $\wt{T}_{Y\sp\sigma}\cong \wt{T}\sprime$.
\end{proposition}
\section{Applications}
In this section,
we consider only lattices without orientation.
For a lattice $T$,
let  $g(T)$ denote the number of isomorphism classes of lattices in
the genus of $T$.
By  Theorem~\ref{thm:main2} and Propositions~\ref{prop:A},~\ref{prop:K3},
we obtain the following:
\begin{corollary}\label{thm:gTA}
Let $X$ be a complex singular abelian surface or a complex singular $K3$ surface,
and $D\subset X$ a reduced effective divisor 
such that the classes of irreducible components of $D$
span $S_X\otimes\Q$.
Then the set of the homeomorphism types of
$(X\sm D)\sp\sigma$ $(\sigma\in \Emb(\C))$
contains at least $g(T_X)$ distinct elements.
\end{corollary}
Another  application is as follows.
By a \emph{plane curve},
we mean a complex reduced (possibly reducible) projective plane curve.
\begin{definition}(\cite{MR1257321}, \cite{MR1980995}, \cite{MR2247887})\label{def:AZP}
A pair $(C, C\sprime)$ of plane curves is said to be an \emph{arithmetic Zariski pair}
if the following hold:

(i) Let $F$ be a homogeneous polynomial defining $C$.
Then there exists $\sigma\in \Emb (\C)$ such that $C\sprime$
is isomorphic (as a plane curve) to $C\sp\sigma:=\{F\sp\sigma=0\}$.

(ii)
There exist  tubular neighborhoods $\TTT\subset \P^2$ of $C$ and $\TTT\sprime\subset \P^2$ of $C\sprime$
such that $(\TTT, C)$ and $(\TTT\sprime, C\sprime)$ are diffeomorphic.

(iii)
$(\P^2, C)$ and $(\P\sp 2, C\sprime)$ are \emph{not} homeomorphic.
\end{definition}
\begin{remark}
The first example of an arithmetic Zariski pair was
discovered  by Artal-Carmona-Cogolludo (\cite{MR1980995},~\cite{MR2247887}) in 
 degree $12$.
\end{remark}
\begin{definition}(\cite{MR805337})
A plane curve $C$ of degree $6$ is called a \emph{maximizing sextic} if 
$C$ has only simple singularities and the  total Milnor number of $C$ 
attains the possible maximum  $19$.
\end{definition}
Let $C$ be a maximizing sextic.
Then,
for any $\sigma\in \Emb(\C)$,
the conjugate plane curve
$C\sp\sigma$ is also a maximizing sextic,
and $(C, C\sp\sigma)$ satisfies the condition (ii) in Definition~\ref{def:AZP},
because simple singularities have no moduli.
We denote by $W_C\to \P^2$  the double covering branched along $C$,
and by $Y_C\to W_C$ the minimal resolution.
Then $Y_C$ is a complex singular $K3$ surface.
We denote by $T[C]$ the transcendental lattice  of $Y_C$.
Let $D_C\subset Y_C$ be the total inverse image
 of $C$ by $Y_C\to W_C\to \P^2$, and we put $U_C:=Y_C\sm D_C$.
Since  the classes of  irreducible components of $D_C$ span $S_{Y_C}\otimes \Q$,
the topological invariant $(B_{U_C}, \beta_{U_C})$ of $U_C$ is isomorphic to $T[C]$.
In particular, 
if 
 $(\P^2, C)$ and $(\P^2, C\sp\sigma)$ are homeomorphic,
then $U_C$ and $U_{C\sp\sigma}$ are also homeomorphic,
and hence 
 $T[C]$ and $T[C\sp\sigma]$ are  isomorphic.
On the other hand,  the set of isomorphism classes of $T[C\sp\sigma] (\sigma \in \Emb(\C))$ 
coincides with the genus of $T[C]$.
 Hence, if $g(T[C])>1$, then there exists $\sigma\in \Emb(\C)$
 such that $(C, C\sp\sigma)$
 is an arithmetic Zariski pair.
 %
%
%
%
Using Yang's algorithm~\cite{MR1387816}, we obtain the following theorem by computer-aided calculation:
{\footnotesize 
\begin{table}
$$
\renewcommand{\arraystretch}{1}
\begin{array}{c|l|l l}
\textrm{No.} & R & \rlap{\textrm{$T[C]$ \;\;and \;\;$T[C\sprime]$}}\mystrutd{3pt}&\\
\hline
1 & E_8 + A_{10}+A_{1} & L[6, 2, 8], & L[2,0,22] \mystruth{10pt}\\
2 & E_8 + A_6 + A_4 +A_1 & L[8, 2, 18], & L[2, 0, 70] \\
3 & E_6 +D_5+A_6+A_2 & L[12, 0, 42], & L[6, 0, 84] \\
4 & E_6+A_{10}+A_3 & L[12,0,22], & L[4,0,66] \\
5 & E_6+A_{10}+A_{2}+A_1 & L[18,6,24], & L[6,0,66]\\
6 & E_6 +A_7 +A_4 + A_2 & L[24, 0, 30], & L[6, 0, 120]\\
7 & E_6+A_6+A_4+A_2+A_1 & L[30, 0, 42], & L[18, 6, 72] \\
8 & D_8+A_{10}+A_{1} & L[6,2,8], & L[2,0,22] \\
9 & D_8+A_6+A_4+A_1 & L[8,2,18], & L[2,0,70]\\
10 & D_7+A_{12} & L[6, 2, 18], & L[2, 0, 52]\\
11 & D_7 + A_8 +A_4 & L[18, 0, 20], & L[2, 0, 180] \\
12 & D_{5}+A_{10}+A_{4} & L[20,0,22], & L[12, 4, 38]\\
13 & D_5 +A_6 +A_5+A_2+A_1 & L[12, 0, 42], & L[6, 0, 84] \\
14 & D_5+A_6+2A_4 & L[20,0,70], & L[10, 0, 140] \\
15 & A_{18}+A_{1} & L[8, 2, 10], & L[2, 0, 38] \\
16 & A_{16}+A_{3} & L[4,0,34], & L[2,0,68] \\
17 & A_{16}+A_{2}+A_{1} & L[10,4,22], & L[6, 0, 34] \\
18 & A_{13}+A_{4}+2A_{1} & L[8,2,18], & L[2,0,70]\\
19 & A_{12}+A_{6}+A_1 & L[8,2,46], & L[2,0,182] \\
20 & A_{12}+A_5+2A_1 & L[12,6,16], & L[4,2,40] \\
21 & A_{12}+A_{4}+A_{2}+A_{1} & L[24, 6, 34], & L[6,0,130] \\
22 & A_{10}+A_{9} & L[10,0,22], & L[2, 0,110] \\
23 & A_{10}+A_{9} & L[8,3,8], & L[2, 1, 28] \\
24 & A_{10}+A_8+A_1 & L[18,0,22], & L[10,2,40] \\
25 & A_{10}+A_7+A_2 & L[22,0,24], & L[6,0,88] \\
26 & A_{10}+A_7+2A_{1} & L[10,2,18], & L[2,0,88] \\
27 & A_{10}+A_6+A_2+A_1 & L[22,0,42], & L[16, 2, 58] \\
28 & A_{10}+A_5+A_3+A_1 & L[12, 0, 22], & L[4,0,66]\\
29 & A_{10}+2A_4+A_1 & L[30, 10, 40], & L[10,0,110] \\
30 & A_{10}+A_4+2A_2+A_1 & L[30, 0, 66], & L[6, 0, 330] \\
31 & A_8+A_6+A_4+A_1 & L[22,4,58], & L[18, 0, 70] \\
32 & A_7+A_6+A_4+A_2 & L[24, 0, 70], & L[6, 0, 280] \\
33 & A_7+A_6+A_4+2A_1 & L[18, 4,32], & L[2, 0, 280]\\
34 & A_7+A_5+A_4+A_2+A_1 & L[24, 0, 30], & L[6, 0, 120] 
\end{array}
$$
\vskip .2cm
\caption{Examples of arithmetic Zariski pairs of maximizing sextics}\label{table:AZP}
\end{table}
}
\begin{theorem}\label{thm:AZP}
There exist  arithmetic Zariski pairs $(C, C\sprime)$ of maximizing sextics
with simple singularities of Dynkin type $R$
for each $R$ in Table~\ref{table:AZP}.
\end{theorem}
The lattices $T[C]$ and $T[C\sprime]$ are also presented
in Table~\ref{table:AZP}.
We denote by $L[2a,b,2c]$ the lattice of rank $2$ represented by the matrix $Q[a,b,c]$.
\begin{remark}
In the previous paper~\cite{AZP},
we have obtained a part of Theorem~\ref{thm:AZP}
by heavily  using results of  Artal-Carmona-Cogolludo~\cite{MR1900779}
and Degtyarev~\cite{degtyarev-2005}.
A detailed account of the algorithm for Theorem~\ref{thm:AZP}
is also given in~\cite{AZP}.
Table~\ref{table:AZP} has been obtained during the calculation
for the results in~\cite{normalK3}.
\end{remark}
\bibliographystyle{plain}

\begin{thebibliography}{10}

\bibitem{MR0349679}
H.~Abelson.
\newblock Topologically distinct conjugate varieties with finite fundamental
  group.
\newblock {\em Topology}, 13:161--176, 1974.

\bibitem{MR1257321}
E.~Artal-Bartolo.
\newblock Sur les couples de {Z}ariski.
\newblock {\em J. Algebraic Geom.}, 3(2):223--247, 1994.

\bibitem{MR1900779}
E.~Artal~Bartolo, J.~Carmona~Ruber, and J.-I.~Cogolludo~Agust{\'{\i}}n.
\newblock On sextic curves with big {M}ilnor number.
\newblock In {\em Trends in singularities}, Trends Math., pages 1--29.
  Birkh\"auser, Basel, 2002.

\bibitem{MR1980995}
E.~Artal~Bartolo, J.~Carmona~Ruber, and J.-I.~Cogolludo~Agust{\'{\i}}n.
\newblock Braid monodromy and topology of plane curves.
\newblock {\em Duke Math. J.}, 118(2):261--278, 2003.


\bibitem{MR2247887}
E.~Artal~Bartolo, J.~Carmona~Ruber, and J.-I.~Cogolludo~Agust{\'{\i}}n.
\newblock Effective invariants of braid monodromy.
\newblock {\em Trans. Amer. Math. Soc.}, 359(1):165--183 (electronic), 2007.

\bibitem{MR0145560}
M.~F. Atiyah and F.~Hirzebruch.
\newblock Analytic cycles on complex manifolds.
\newblock {\em Topology}, 1:25--45, 1962.


\bibitem{MR1028322}
D.~A. Cox.
\newblock {\em Primes of the form {$x\sp 2 + ny\sp 2$}}.
\newblock A Wiley-Interscience Publication. John Wiley \& Sons Inc., New York,
  1989.

\bibitem{degtyarev-2005}
A.~Degtyarev.
\newblock On deformations of singular plane sextics.
\newblock Preprint,  2005. To appear in J. Algebraic Geom.
\newblock  http://arxiv.org/abs/math.AG/0511379.


\bibitem{MR1644323}
W.~Fulton.
\newblock {\em Intersection theory}, volume~2 of {\em Ergebnisse der Mathematik
  und ihrer Grenzgebiete. 3. Folge. A Series of Modern Surveys in Mathematics.}
\newblock Springer-Verlag, Berlin, second edition, 1998.

\bibitem{MR890960}
S.~Lang.
\newblock {\em Elliptic functions}, volume 112 of {\em Graduate Texts in
  Mathematics}.
\newblock Springer-Verlag, New York, second edition, 1987.

\bibitem{MR525944}
V.~V. Nikulin.
\newblock Integer symmetric bilinear forms and some of their geometric
  applications.
\newblock {\em Izv. Akad. Nauk SSSR Ser. Mat.}, 43(1):111--177, 238, 1979.
\newblock English translation: Math USSR-Izv. 14 (1979), no. 1, 103--167
  (1980).

\bibitem{MR805337}
U.~Persson.
\newblock Double sextics and singular {$K$}-{$3$} surfaces.
\newblock In {\em Algebraic geometry, Sitges (Barcelona), 1983}, volume 1124 of
  {\em Lecture Notes in Math.}, pages 262--328. Springer, Berlin, 1985.

\bibitem{MS}
M.~Sch\"utt.
\newblock Fields of definition for singular {$K3$} surfaces. 
\newblock Preprint, 2006.
To appear in Communications in Number Theory and Physics.
\newblock http://arxiv.org/abs/math.AG/0612396.

\bibitem{MR0166197}
J.-P.~Serre.
\newblock Exemples de vari\'et\'es projectives conjugu\'ees non hom\'eomorphes.
\newblock {\em C. R. Acad. Sci. Paris}, 258:4194--4196, 1964.


\bibitem{normalK3}
I.~Shimada.
\newblock On normal {$K3$} surfaces.
Preprint, 2006. To appear in Michigan Math.~J.
\newblock http://arxiv.org/abs/math.AG/0607450.


\bibitem{AZP}
I.~Shimada.
\newblock On arithmetic {Z}ariski pairs in degree $6$.
\newblock Preprint, 2006. To appear in Adv. Geom.  
\newblock http://arxiv.org/abs/math.AG/0611596.


\bibitem{tssK3}
I.~Shimada.
\newblock Transcendental lattices and supersingular reduction lattices of a
  singular {$K3$} surface. 
\newblock Preprint,  2006. To appear in Trans. Amer. Math. Soc.\\
\newblock  http://arxiv.org/abs/math.AG/0611208.


\bibitem{MR0314766}
G.~Shimura.
\newblock {\em Introduction to the arithmetic theory of automorphic functions}.
\newblock Publications of the Mathematical Society of Japan, No. 11. Iwanami
  Shoten, Publishers, Tokyo, 1971.

\bibitem{MR0441982}
T.~Shioda and H.~Inose.
\newblock On singular {$K3$} surfaces.
\newblock In {\em Complex analysis and algebraic geometry}, pages 119--136.
  Iwanami Shoten, Tokyo, 1977.

\bibitem{MR0382289}
T.~Shioda and N.~Mitani.
\newblock Singular abelian surfaces and binary quadratic forms.
\newblock In {\em Classification of algebraic varieties and compact complex
  manifolds}, pages 259--287. Lecture Notes in Math., Vol. 412. Springer,
  Berlin, 1974.

\bibitem{MR1423033}
B.~Totaro.
\newblock Torsion algebraic cycles and complex cobordism.
\newblock {\em J. Amer. Math. Soc.}, 10(2):467--493, 1997.

\bibitem{MR1967689}
C.~Voisin.
\newblock {\em Hodge theory and complex algebraic geometry. {I}}, volume~76 of
  {\em Cambridge Studies in Advanced Mathematics}.
\newblock Cambridge University Press, Cambridge, 2002.

\bibitem{MR1387816}
Jin-Gen Yang.
\newblock Sextic curves with simple singularities.
\newblock {\em Tohoku Math. J. (2)}, 48(2):203--227, 1996.


\end{thebibliography}

\def\cprime{$'$} \def\cprime{$'$} \def\cprime{$'$}

\end{document}